\documentclass[twoside,11pt]{article}
\usepackage{pst-all}
\usepackage{pstricks}
\usepackage{url}
\usepackage[pagebackref=true,colorlinks=false,urlcolor=blue]{hyperref}
\usepackage{multirow}
\usepackage{natbib}
\usepackage{graphics}
\usepackage{amsfonts}
\usepackage{amsmath}
\usepackage{float}
\usepackage{pgf}
\newcommand{\email}[1]{\href{mailto:#1}{#1}}

\title{Pickup and Delivery Problem with Transfers}

\author{Santiago Hincapie Potes \thanks{Mathematical Engineering Student,
    \email{shinca12@eafit.edu.co}, Universidad EAFIT,\newline
    Medellín, Colombia.}, 
Catalina Lesmes Ramirez\thanks{Mathematical Engineering Student,
\email{clesmes@eafit.edu.co}, Universidad EAFIT,\newline Medellín,
Colombia.} \hspace{-11pt}}

\date{}

\begin{document}
\maketitle 

\begin{abstract}
  In this article we will be talking about Pickup and Delivery Problems with
  Transfers (PDP-T), the main idea is to show first the Pickup and Delivery
  Problem (PDP), which is the problem from which the PDP-T derives from, later
  show a description of the problem and the mathematical formulation to get a
  clear view of what we are going to solve, which leads us to present different
  ways of solving this kind of problems. 
\end{abstract}


\section{Introduction}\label{sec_intro}

In many logistics problems we require to pick some products in one place and
take them to another place, which is why we define the pickup and
delivery problem (PDP). In this problem a set of vehicles pick up and deliver
a set of items. The goal is to deliver the items at the lowest
cost while obeying a set of constraints, such as time windows and capacities.
The PDP is a well-studied, NP-hard problem, so approximation algorithms and
heuristics have been developed to address variants of the PDP.
There are many techniques we can use to solve the PDP problem, such as
genetic algorithms, various metaheuristics, taboo search heuristics and
branch and cut algorithms.

To the PDP certain transfers can be added, which means that the vehicle which is delivering
the product transfers it to another vehicle before the delivery is done. By adding those transfers we have the PDP with Transfers (PDP-T), in which we consider
transferring items between vehicles. We can convert the PDP problem into a
PDP-T problem by adding some variables and constraints.

This article is organized as follows: In Section 1 an  introduction is presented which gives a general notion of the problem, later in Section 2 we give a brief literature review were some examples of solution methods for the PDP-T are presented, after in Section 3 we give the problem description, were all the variables, parameters and constraints of the mathematical model are explained in detail, next in Section 4 the solution algorithms we give to the problem are explained, later in Section 5 we present the results of the computational experimental and finally in Sections 6 and 7 we give some conclusions and future work suggestions. 
\section{Literature Review}\label{sec_SoA}

As we saw before, there are different methods to solve the PDP-T.

As an example a branch and cut method using Benders Decomposition can be used. In this method,
the set of constraints is decomposed into pure integer and mixed constraints,
and then a branch-and-cut procedure is applied to the resulting pure integer
problem, by using real variables and constraints related as cut generators.
The key on the success of this method is that those constraints defined by a
logical sentence are not modeled using the big-M technique, as usual in a
branch-and-bound methodology behind the original PDP formulation. This method
may be applied only when the objective function is either pure real or pure
integer \citep{Cortes2010}.

Another method which solves the problem is a Very Large Neighborhood Search with Transfers
(VLNS-T) based on the Adaptive Very Large Neighborhood Search (VLNS).
The VLNS algorithm  uses  simulated  annealing  to  randomly  choose
neighboring schedules and iteratively improve the schedule. Neighboring
schedules are formed by removing random items and reinserting them with
heuristics. 
VLNS-T is based on the VLNS algorithm for the PDP without transfers, a
variant of simulated annealing in which the neighborhood of states is "very
large". In this case we remove random items from the schedule and then reinsert
them with multiple heuristics to find neighbors \citep{Cata}.

There are also some special cases of the PDP-T, as the Pickup and
Delivery Problem with Shuttle Routes (PDP-S) which relies on a structured
network with two categories of routes. Pickup routes visit a set of pickup
points independently of their delivery points and end at one delivery point.
Shuttle routes are direct trips between two delivery points. Requests can be
transported in one leg (pickup route) or two legs (pickup route plus shuttle
route) to their delivery point. The PDP-S applies to transportation systems
with a multitude of pickup points and a few, common delivery points
\citep{Masson2010}.
\section{Problem Description and Mathematic Formulation}\label{sec_mathmod}
In the PDP, a heterogeneous vehicle fleet based at multiple terminals must
satisfy a set of transportation requests. Each request is defined by a pickup
point, a corresponding delivery point, and a demand to be transported between
these locations \citep{Parragh2008}. Now, PDP-T allows the option for passengers
to transfer between vehicles, provided that the locations of the transfer
points are fixed and known.\\
\vspace{-8mm}
\subsection{Problem Formulation}
Let $N$ be the set of transportation requests. For each transportation request
$i\in N$ a load of size $\bar{q_i}\in \mathbb{N}$ has to be transported from a
set of origins $N^+_i$ to a set of destinations $N_i^-$. \\
Each load is divided as follows
$\displaystyle\bar{q_i} = \sum_{j\in N_i^+}q_j = - \sum_{j\in N_i^-} q_j$,
i.e., positive quantities for pickups and negative quantities for deliveries.\\
Define $\displaystyle N^+ := \cup_{i\in N} N_i^+$ as the set of all origins and
$\displaystyle N^- := \cup_{i\in N} N_i^-$ as the set of all destinations.
Let $V := N^+\cup N^-$.\\
Furthermore, let $M$ be the set of vehicles. Each vehicle $k\in M$ has a
capacity $Q_k \in \mathbb{N}$ a start location $k^+$ and an end location $k^-$.
\\Define $M^+ := \{k^+|k\in M\}$ as the set of start locations and
$M^- := \{k^-|k\in M\}$ as the set of end locations. Let $W := M^+ \cup M^-$.\\
For all $i, j \in V \cup  W$ let $d_{ij}$ denote the travel distance, $t_{ij}$
the travel time and $c_{ij}$ the travel cost. Note that the dwell time at
origins and destinations can be easily incorporated in the travel time and
therefore will not be considered explicitly.\\

To formulate the PDP as a mathematical program we introduce four types of
variables:
for $i\in N$ and $k\in M$:
$$
z_i^k =
\begin{cases}
  1 & \text{if transportation request } i \text{ is assigned to vehicle } k \\
  0 & \text{In other case}
\end{cases}
$$
for $(i, j)\in (V \times V)\cup \{(k^+,j)|j\in V\} \cup \{(j,k^-)|j\in V\}$
and $k\in M$:
$$
x_{ij}^k = \begin{cases} 1 & \text{if vehicle }k\  \text{travels from location }
i\  \text{to location } j \\ 0 & \text{In other case}  \end{cases}
$$
$D_i$  with $(i \in  V\cup W)$, specifying the departure time at vertex $i$ and
$y_i$ with $(i \in  V\cup W)$, specifying the load of the vehicle arriving at
vertex $i$
\citep{Parragh2008}.\\
All this information can be summarized in Table 1:
\begin{table}[H]
\caption{Parameters and variables found in 
\citep{Parragh2008} for the model}
\begin{center}
\begin{tabular}{ll}
  \hline
  Object   & Meaning                                             \\
  \hline
  $M$      & Set of vehicles                                     \\
  $C$      & Set of requests                                     \\
  $T$      & Set of transference point                           \\
  $M^+$    & Set of origin depots for vehicles                   \\
  $M^-$    & Set of destination depots for vehicles              \\
  $N^+$    & Set of origin nodes for requests                    \\
  $N^-$    & Set of destination nodes for requests               \\
  $N$      & Set of nodes associated with requests               \\
  $V$      & Set of nodes                                        \\
  $q_{ij}$ & Size of request $i \in C$                           \\
  $Q_k$    & Capacity of vehicle $k \in K$                       \\
  $t_{ij}$ & Minimum ride time from node $i$ to node $j$         \\
  $c_{ij}$ & The travel cost                                     \\
  $d_{ij}$ & The travel distance                                 \\
  $z_i^k$  & bind transportation request and vehicles            \\
  $x_{ij}$ & bind routes and vehicles information                \\
  $D_i$    & specifying the departure time at specific vertex    \\
  $y_i$    & specifying the load of the vehicle arriving         \\
  \hline
\end{tabular}
\end{center}
\end{table}

Now the model is:
\begin{align}
& {\text{minimize}}
& & \sum_{i,j\in V \cup W} dij && \\
& \text{subject to}
& & \sum_{k\in M} z_i^k = 1 & \text{for all }& i\in N \\
&&& \sum_{j\in V\cup W} x_{il}^k = z_{i}^k & \text{for all } &
    i \in N,l\in N^+_i\cup N_i^- k\in M\\
&&& \sum_{j\in V\cup\{k^-\}} x^k_{k^+ j} = 1 & \text{for all }& k\in M\\
&&& \sum_{j\in V\cup\{k^+\}} x^k_{k^- j} = 1 & \text{for all }& k\in M\\
&&& D_{k^+} = 0 & \text{for all }& k\in M\\
&&& D_p \leq D_q & \text{for all }& i\in N, p\in N_i^+ q\in N_i^-\\
&&& x_{ij}^k = 1\Rightarrow D_i + t_{ij} \leq D_j & \text{for all } &
    i,j \in V\cup W, k\in M\\
&&& y_{k^+} = 0 & \text{for all }& k\in M\\
&&& y_l \leq \sum_{k\in M} Q_kz_i^k & \text{for all } &
    i\in N, l\in N_i^+\cup N_i^-\\
&&& D_i \geq 0 & \text{for all }& i\in V\cup W\\
&&& y_i \geq 0 & \text{for all }& i\in V\cup W
\end{align}

Constraint (2) ensures that each transportation request is assigned to exactly
one vehicle. By constraint (3) a vehicle only enters or leaves a location $l$
if it is an origin or a destination of a transportation request assigned to that
vehicle. The next (4) and (5) make sure that each vehicle starts and ends at
the correct place. Also the constraints number (6), (7), (8) and (11) together
form the precedence constraints the others together form the capacity
constraints.  Constraints (9), (10) and (12) together form the capacity
constraints.\\
This mathematical program, models the PDP. Now for a PDP-T model in the
literature it can bee seen that the transfer point is introduced, and the procedure for the extended
model is to iterative add constrains that involve this transfer points. The process of inserting the transfers will be described in Section 4.

\section{Solution Algorithms}
In this section we present solution algorithm for the PDP-T.
In general, for this part, we first solve the PDP using different methods (which would be described in this section) and then we add the transship option to the PDP solution, transship option always increases the computational time but makes the solution better. In subsection 4.1 we present two constructive methods: a greedy approach and multistart heuristic in \cite{Takoudjou2012}, then in subsection 4.2 we present a random method GRASP which is applied to the constructive methods mentioned before, later in subsection 4.3 we present three local search methods: Variable neighborhood descent, Adaptive large neighborhood search and simulated annealing, finally in subsection 4.4 we present two hybrid methods which combine local search methods and genetic algorithms. 
\subsection{Constructive Methods}
\subsubsection{Greedy approach}
We implemented the greedy approach described in \cite{coltin2014multi}, which basically iterates through every item and vehicle and then inserts the best pickup and delivery action into the schedule, it always chooses the option that increases the least the cost, we repeat the process until no unassigned items remain. Then we add the transfer option with another greedy idea, inspired on Clarke saving algorithm, this method computes the savings of the transfer, and it improves the solution, then we modify the solution by adding this new note.
\subsubsection*{ Clarke and Wright Algorithm}
The Clarke and Wright savings algorithm is one of the most known heuristic for VRP. It applies to problems for which the number of vehicles is not fixed (it is a decision variable), and it works equally well for both directed and undirected problems. When two routes $(0,...,i,0)$ and $(0,j,...,0)$ can feasibly be merged into a single route $(0,...,i,j,...,0)$, a distance saving $s_{ij} = c_{i0} + c_{0j} - c_{ij}$ is generated. The algorithm works as it follows:\\

\textbf{Step 1. Savings computation} 
\begin{itemize}
\item Compute the savings $s_{ij} = c_{i0} + c_{0j} - c_{ij}$ for $i,j = 1,...,n$ and $i \neq j$.
\item Create n vehicle routes $(0,i,0)$ for $i = 1,...,n$.
\item Order the savings in a non increasing fashion.
\end{itemize} 

\textbf{Step 2. Route Extension (Sequential version)}
\begin{itemize}
\item Consider in turn each route $(0,i,...,j,0)$.
\item Determine the first saving $s_{ki}$ or $s_{jl}$ that can feasibly be used to merge the current route with another route ending with $(k,0)$ or starting with $(0,l)$.
\item Implement the merge and repeat this operation to the current route. 
\item If not feasible merge exists, consider the next route and reapply the same operations.
\item Stop when not route merge is feasible.
\end{itemize}

\subsubsection{Hybrid Multistart Heuristic}
The PDPT is a hard combinatorial optimization
problem. The heuristics for this problem
must avoid being trapped by local optimum. To overcome
local optimality, a diversification procedure is
needed. So we use the transshipment heuristic from the hybrid heuristic for the PDP found in \cite{Takoudjou2012} to obtain a solution to the PDP-T.\\
At each iteration, the transshipment heuristic is used
to improve the PDP solution and obtain a solution
to the PDP-T. From the current PDP solution, each
request ($i$, $i + n$) $\in$ R; is removed from the solution.
Then, ($i$, $i + n$) is split into two different requests
($i$, $e_{t}$) and ($s_{t}$, $i + n$), where $e_{t}$ and $s_{t}$ are the inbound/outbound doors of a transshipment point
$t$ $\in$ T. The best reinsertion cost of ($i$, $i + n$) in the
solution is computed. The best insertion cost to insert
($i$, $e_{t}$) following by insertion of ($s_{t}$, $i + n$) in the
solution is computed. The cost to insert ($s_{t}$, $i+n$) following
by the insertion of ($i$, $e_{t}$) at their best position
is computed. Between the three possibilities, the insertion
or the reinsertions offering the minimum cost
is performed.
As described in \cite{Takoudjou2012}, this hybrid approach merges many ideas. The
solving principle is based on the following steps:
\begin{itemize}
  \item \textbf{An initial PDP solution (Sol) is calculated by gradually
                inserting the requests in routes associated to vehicles.}\\
                The initial solution is consisting of a single route containing
                only the starting point and the ending point of each route.
                Requests are then successively introduced into the tour of a
                vehicle offering the minimal increase of the cost of transport. \\
  \item \textbf{Transshipment is then used to “destroy” and “repair” PDP
                solution to obtain a better solution (PDPT solution).}\\
                At each iteration, the transshipment heuristic is used to improve
                the PDP solution and obtain a solution to the PDPT. From the
                current PDP solution, each request; is removed from the solution.
                Then, the request is split into two different requests, the first
                go from the start point to the best transfer point, the other go
                from the transfer point to the ending point.
\end{itemize}
\subsection{Random Search Method}
\subsubsection{GRASP}
GRASP is a multi-start metaheuristic for combinatorial optimization problems, in which each iteration consists of two phases: construction and local search. The construction phase builds a feasible solution, whose neighborhood is investigated until a local minimum is found during the local search phase. The best overall solution is kept as the result. Here we will only use one phase, which is the construction of the solution \cite{grasp}.

\subsection{Local Search Methods}
\subsubsection{Variable Neighborhood Descent}
The VND is used
to explore the neighborhood of the current solution,
which is based on three operators defined below: ADR, RNR and SWR. With these
three operators, we can explore the solution space
more intensively. An operator is a move that transforms
one solution to another with small modifications. In best improvement operators the
whole neighborhood is analyzed and the best solution
is kept. In first improvement local search the first
better solution found in the neighborhood is kept \citep{Takoudjou2012}.
\paragraph{SWR}(Swap requests between routes)
In SWR neighborhood, two requests belonging to two
different routes are exchanged together provided
that all PDP constraints are satisfied. 
\paragraph{RNR}(Remove and insert a request)
In RNR neighborhood,
provided that all PDP constraints are
satisfied, request belonging to one route is removed
and inserted in another route. 
\paragraph{ADR} (Advance or delay a request)
In ADR neighborhood,
requests belonging to one given route
is advanced or delayed in the same route if PDP
constraints are all satisfied.

\subsubsection{Adaptive Large Neighborhood Search}
The general idea is to remove some requests
from routes in the current solution and then reinsert them elsewhere
to arrive at a new solution. At each iteration, one each of
several removal and insertion heuristics is randomly selected based
on their designated weights. Letting the set of removal heuristics be
${rh_{1},..., rh_{nr}}$ and their corresponding weights be ${rw_{1},..., rw_{nr}}$, $rh_{i}$
is chosen with probability $\frac{rw_{i}}{\sum_{j=1}^{nr}rw_{j}}$, where $nr$ is the number of options. Similarly for the insertion heuristics ${ih_{1},..., ih_{ni}}$ and their
corresponding weights ${iw_{1},..., iw_{ni}}$, the probability of $ih_{i}$ being $\frac{ih_{i}}{\sum_{j=1}^{n}ih_{j}}$
selected is, where $ni$ is the number of options. After
each iteration, if the solution realized is feasible and better than the
incumbent, the latter is updated. In any case, the process is repeated
until a predefined number of iterations $n_{II}^{max}$ is reached \citep{Masson2012}.

\subsubsection{Simulated Annealing}
Simulated annealing is a metaheuristic that begins at some state, and chooses a random “neighbor” of that state. With probability accept $(e, e_{0}, t)$ the new state is accepted as the current state, where $e$ is the “energy” of the current state, $e_{0}$ is the energy of the new state, and $t$ is the temperature, or the fraction of iterations of the algorithm currently completed. If the new state is rejected we remain at the current state and repeat with a new neighbor. The algorithm continues either for a
fixed number of iterations or until the energy crosses some threshold, when the best solution that has been encountered thus far is returned \citep{Coltin2014}.

\subsection{Hybrid Methods}
In this subsection we present two hybrid methods, both of this methods are based on the use of a local search methods and genetic algorithms but they possess different hybridizations. The first algorithm uses simulated annealing for the population selection and for its mutation variance of the probability, this is made with the purpose of achieves a better convergence, next to a diverse population. The second hybrid method uses the taboo list to avoid mating of chromosomes that have low Minkowski in order to keep diversity, also use this information in order to change mutation probability. We show in precise detail how the genetic algorithm was made, so in subsection 4.4.1 we show the encoding of the problem, then we describe the initial population in subsection 4.4.2, which we got in two different ways, we have a totally randomized initial solution and a constructive solution, later we select the population to be matched which is described in subsection 4.4.3, then in subsection 4.4.4 we describe the crossover between the population, after that, in subsection 4.4.5 we describe the mutation for both first and second methods, then in subsection 4.4.6 we explain how the new population is chosen for the first method and the second method, and finally in subsection 4.4.7 we explain the stop criterion.

\subsubsection{Encoding}
Given the structure of the problem, the solutions can be represented as a set of $N$ paths, one for each vehicle, considering that it is important to consider the transfers as much as the possibility of a node recurrence, it was opted to represent the solution as a directed multi-graph, in which the set of vertexes represent some of the visited nodes and there exists a connexion between two vertexes if one is a successor from the other, the direction goes in the way of the route, additionally each one of the multi-sets is found reserved for each route. To represent each multi-graph, a multidimensional matrix was used, which positions $i, j, k$ possessed the locality distance $i$ to the locality $j$, if $j$ is the successor of $i$ in route $k$.
For the processes of the genetic algorithm to each element of the solution matrix a derivate from the logistic function was applied (fixing the value from $0$ to $0$), meaning that to each element a function as is shown below was applied: 
$$ 
f(x) = 
\begin{cases}
\frac{e^{-x}}{(1 + e^{-x})^2} & \text{if } x > 0 \\
0                             & \text{otherwise}
\end{cases}
$$ 

this function was applied with the purpose of giving more importance to the short routes, in this way the intention is to preserve this property.
The complexity in space of this encoding is $O(N|path|^2)$.

\subsubsection{Initial Population}
For the computational experiments we performers with two district initial population methods, a totally randomized approach and a set of district greedy random solutions.
\subsubsection*{Totally Randomized Initial Solution}
A set of 256 random multi-graphs were generated using the model of a modification of the algorithm from Erdős-Rényi as described in \citep{initpopul}. It was opted to generate random graphs in place of random matrices, due to the fact that the solutions generated this way are usually more feasible with a higher probability, additionally they use to be better.
\subsubsection*{Constructive Solution}
We construct a set of 256 district random greedy solution using the algorithm described in the section 4.2.

\subsubsection{Selection}
For the selection  part we use Baker’s stochastic universal selection (SUS) \citep{SUS}, this algorithm uses a single random value to sample all the solutions by choosing them at evenly spaced intervals. This gives weaker members of the population (according to their fitness) a chance to be chosen and thus reduces the unfair nature of fitness-proportional selection methods, also the experimental work by Hancock \citep{Hancock1994} clearly demonstrates the superiority of this approach.

\subsubsection*{Fitness Function}
For a $K$ matrix its fitness function is the sum of the total distance of each one of its routes, in case that the solution is not feasible (it does not satisfy all the request) a $\infty$ value is assigned.

\subsubsection{Crossover}
To make the crossover the solutions which were previously selected were matched and a new matrix was generated in which its entries were uniform random numbers. The resultant son matrix  possesses the gen of its first parent in the position $i, j, k$ if the generated random matrix possesses a lower number to $0.5$ in the corresponding position, in the other cases it possesses the gen of its other parent.

\subsubsection{Mutation}
\paragraph{First method}
For mutation operator we use a shift operator, that shifts the value of two elements of the matrix depending on a parameter that we varied according to the diversity in the population (measured in terms of the fitness variance) and actual solution temperature (describe in the next step).
\paragraph{Second Method}
For mutation operator we also use a shift operator but in this method we uses a fixed probability $\alpha = 0.05$

\subsubsection{New population}
\paragraph{First method}
For the actualization of the population, we perform a probabilistic replacement test, in this approach the solutions are accepted into the new generation using a Boltzmann type calculation. At each generation the temperature is decreased, thereby decreasing the probability of accepting higher energy individuals.
\paragraph{Second method}
For the population update, we compute the Minkowski distance between population members and uses taboo list to avoid nears solutions, in order to keep population diversity.

\subsubsection{Stop criterion}
For the stop criterion we compared the variance between the solutions of the last $4$ generations, this with the purpose of ending the algorithm when it reaches an stationary state.
\section{Computational Experimentation}\label{sec_results}
In this section we show the computational results given by the algorithms described in section 4 that solve the PDP-T. In subsection 5.1, we present tables with the instance size, the number of vehicles, optimal solutions and the computational time, then we compare them to see which is the best heuristic used.

All algorithms were coded on python 3.5.3 running on Intel Xeon E5-1620V4
, 64-bit processor with 32 gigabyte RAM and Ubuntu 17.10. We used problem instances
from those in for the PDP, which are related to the well-known Solomon instances.
The datasets are available at\\http://www.sintef.no/Projectweb/TOP/PDPTW/ . For
each instance we introduce a random number of transfer points, in a random node,
to become the original PDP instance in a PDP-T instance.

Each data set contains $X - Y$ hence, the lengths of the arcs in the network are
$L_2$ distances. For the locations of the vehicle depots, the 
origin and final depots for each vehicle were randomly generated, scattered over the region formed by
the nodes. The number of vehicles are also random generated, and finally we
ignore the time windows restriction because are not consider in the model.
\subsection{Computational results}\label{sec_eval}
In this subsection we show the results of all the heuristics mentioned before and compare them. 
\subsubsection{Constructive Methods}
In this subsection we show the tables of the computational results for the constructive methods, in Table 2 we present the results for the Greedy approach and in Table 3 we show the results for the MULTISTART method.
\begin{table}[H]
\caption{Greedy approach}
\begin{tabular}{lllll}
  \hline
  Instance   & Instance Size & Number of vehicles & Optimal Solution & Time      \\
  \hline
  lc204      & 100           & 3                  & 2001.31       & 28.23     \\
  lc2\_2\_3  & 200           & 6                  & 4920.43       & 50.48     \\
  lc2\_4\_5s & 400           & 12                 & 7235.92       & 123.43    \\
  LC2\_8\_5  & 800           & 25                 & 226748.20     & 240.12    \\
  \hline
\end{tabular}
\end{table}

\begin{table}[H]
\caption{MULTISTART}
\begin{tabular}{lllll}
  \hline
  Instance & Instance Size & Number of vehicles & Optimal Solution & Time        \\
  \hline
  lc204      & 100           & 3                  & 1825.83       & 40.81     \\
  lc2\_2\_3  & 200           & 6                  & 4378.92       & 76.21     \\
  lc2\_4\_5s & 400           & 12                 & 6891.06       & 173.79    \\
  LC2\_8\_5  & 800           & 25                 & 201234.73     & 301.83    \\
  \hline
\end{tabular}
\end{table}
As we can see, the MULTISTART method gives better optimal solutions for the PDP-T, but at the same time it requires more computational time.
\subsubsection{Random Search Method}
In this subsection the results for the GRASP which is the random search method are shown, as we mentioned in Section 4, we add GRASP to the constructive methods. In Table 4, the greedy approach with GRASP is shown and in Table 5 we have the results for the MULTISTART with GRASP.
\begin{table}[H]
\caption{Greedy approach (GRASP)}
\begin{tabular}{lllll}
  \hline
  Instance & Instance Size & Number of vehicles & Optimal Solution & Time        \\
  \hline
  lc204      & 100           & 3                  & 1713.01       & 54.12     \\
  lc2\_2\_3  & 200           & 6                  & 4001.12       & 80.15     \\
  lc2\_4\_5s & 400           & 12                 & 6481.10       & 196.92    \\
  LC2\_8\_5  & 800           & 25                 & 198917.87     & 320.16    \\
  \hline
\end{tabular}
\end{table}

\begin{table}[H]
\caption{MULTISTART (GRASP)}
\begin{tabular}{lllll}
  \hline
  Instance & Instance Size & Number of vehicles & Optimal Solution & Time        \\
  \hline
  lc204      & 100           & 3                  & 1521.31       & 58.63     \\
  lc2\_2\_3  & 200           & 6                  & 3832.12       & 89.53     \\
  lc2\_4\_5s & 400           & 12                 & 6312.18       & 200.12    \\
  LC2\_8\_5  & 800           & 25                 & 190143.12     & 354.10    \\
  \hline
\end{tabular}
\end{table}
If we compare Tables 4 and 5 we can see that the optimal solutions are given by the MULTISTART (GRASP), but it also takes more computational time than the Greedy approach (GRASP). This results are coherent, since they follow the same behavior as we mentioned in subsection 5.1.1, but we can also see that when we add GRASP to our constructive methods the optimal solutions are better, even thought they take more computational time.
\subsubsection{Local Search Methods}
In this subsection we show tables comparing all the results, plus one more, the other is called "MIX", this table refers to a union of the VND and ALNS methods, we wanted to see what happened if we merged two of them. So in Table 6 we present the ALNS, in Table 7 the VND, in Table 8 the SA and in table 9 the mix between the VND and the ALNS. Also in this subsection we also wanted to show a 3D plot, showing the results of all the heuristics including the initial solution, this plot is shown in Figure 1.
\begin{table}[H]
\begin{center}
\caption{ALNS}
\begin{tabular}{lllll}
  \hline
  Instance & Instance Size & Number of vehicles & Optimal Solution & Time        \\
  \hline
  lc204      & 100           & 3                  & 1154.21       & 160.32    \\
  lc2\_2\_3  & 200           & 6                  & 3213.31       & 233.21    \\
  lc2\_4\_5s & 400           & 12                 & 5421.32       & 312.23      \\
  LC2\_8\_5  & 800           & 25                 & 182313.32     & 600.12      \\
  \hline
\end{tabular}
\end{center}
\end{table}
\begin{table}[H]
\begin{center}
\caption{VND}
\begin{tabular}{lllll}
  \hline
  Instance & Instance Size & Number of vehicles & Optimal Solution & Time        \\
  \hline
  lc204      & 100           & 3                  & 1032.12       & 187.21    \\
  lc2\_2\_3  & 200           & 6                  & 3102.23       & 254.12    \\
  lc2\_4\_5s & 400           & 12                 & 5250.12       & 352.21    \\
  LC2\_8\_5  & 800           & 25                 & 173213.12     & 721.48    \\
  \hline
\end{tabular}
\end{center}
\end{table}
    
\begin{table}[H]
\begin{center}
\caption{SA}
\begin{tabular}{lllll}
  \hline
  Instance & Instance Size & Number of vehicles & Optimal Solution & Time        \\
  \hline
  lc204      & 100           & 3                  & 1321.23       & 140.32    \\
  lc2\_2\_3  & 200           & 6                  & 3543.34       & 221.34    \\
  lc2\_4\_5s & 400           & 12                 & 5732.12       & 301.31    \\
  LC2\_8\_5  & 800           & 25                 & 182313.32     & 572.12    \\
  \hline
\end{tabular}
\end{center}
\end{table}
    
\begin{table}[H]
\begin{center}
\caption{MIX}
\begin{tabular}{lllll}
  \hline
  Instance & Instance Size & Number of vehicles & Optimal Solution & Time        \\
  \hline
  lc204      & 100           & 3                  & 1023.62       & 350.14    \\
  lc2\_2\_3  & 200           & 6                  & 3102.23       & 500.51    \\
  lc2\_4\_5s & 400           & 12                 & 5291.17       & 698.92    \\
  LC2\_8\_5  & 800           & 25                 & 177843.01     & 1425.92   \\
  \hline
\end{tabular}
\end{center}
\end{table}
\begin{figure}[H]
\centering
\scalebox{0.4}{\input{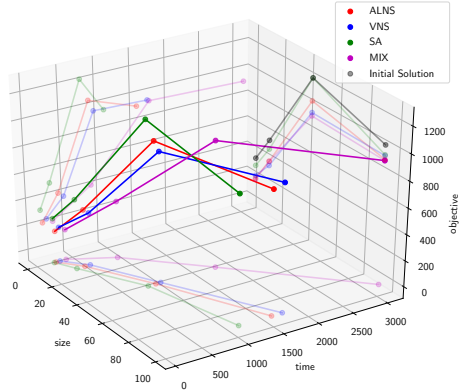}}
\caption{Comparison of the heuristic methods}
\end{figure}

In Figure 1, we can see in the plot the time, the objective function and the size of the instance in contrast between each other, is a comparison of all methods and the initial solution.\\

We can see from Tables 6 through 9 and Figure 1 that the best algorithm as seen in computational time is the Simulated Annealing and that best algorithm as seen in optimal solution is the mix between the ALNS and the VND. As until now, our best results are given with the random search algorithms if we compare them with the local search algorithms we can say that the best results are given by the local search algorithms specifically with the mix of the ALNS and the VNS. For this mix, we can see a very interesting result, even thought it was the best result of all, we should expect it gave a much better result, but the difference from the initial solution is not that wide. 
\subsubsection{Hybrid Methods}
In this subsection we show the results given by the algorithms mentioned in subsection 4.4, it is important to mention that the reaction of the initial population and the associated individual operations such as encoding were computed in parallel. So, in Tables 10 and 11 we show the results for the genetic algorithms and simulated annealing with both the constructive and the random initial solution, and in Tables 12 and 13 we show the results for the genetic algorithms and taboo with both the constructive and the random initial solutions. 

\begin{table}[H]
\begin{center}
\caption{Genetic algorithm \& SA, with random initial solution}
\begin{tabular}{lllll}
  \hline
  Instance & Instance Size & Number of vehicles & Optimal Solution & Time     \\
  \hline
  lc204      & 100           & 3                  & 720.60        & 789.12    \\
  lc2\_2\_3  & 200           & 6                  & 2454.12       & 1325.21   \\
  lc2\_4\_5s & 400           & 12                 & 4843.93       & 2012.32   \\
  LC2\_8\_5  & 800           & 25                 & 14794.27      & 3024.15   \\
  \hline
\end{tabular}
\end{center}
\end{table}

\begin{table}[H]
\begin{center}
\caption{Genetic algorithm \& SA, with Constructive solution}
\begin{tabular}{lllll}
  \hline
  Instance & Instance Size & Number of vehicles & Optimal Solution & Time        \\
  \hline
  lc204      & 100           & 3                  & 673.12	      & 1254.65   \\
  lc2\_2\_3  & 200           & 6                  & 2021.23       & 1743.54   \\
  lc2\_4\_5s & 400           & 12                 & 4654.63       & 2433.33   \\
  LC2\_8\_5  & 800           & 25                 & 13647.12      & 3502.12   \\
  \hline
\end{tabular}
\end{center}
\end{table}

\begin{table}[H]
\begin{center}
\caption{Genetic algorithm \& Taboo, with random initial solution}
\begin{tabular}{lllll}
  \hline
  Instance   & Instance Size & Number of vehicles & Optimal Solution & Time      \\
  \hline
  lc204      & 100           & 3                  & 690.43	      & 823.43    \\
  lc2\_2\_3  & 200           & 6                  & 2144.33       & 1513.43   \\
  lc2\_4\_5s & 400           & 12                 & 4136.38       & 2243.43   \\
  LC2\_8\_5  & 800           & 25                 & 13474.76      & 3324.69   \\
  \hline
\end{tabular}
\end{center}
\end{table}

\begin{table}[H]
\begin{center}
\caption{Genetic algorithm \& Taboo, with Constructive solution}
\begin{tabular}{lllll}
  \hline
  Instance & Instance Size & Number of vehicles & Optimal Solution & Time        \\
  \hline
  lc204      & 100           & 3                  & 610.60	      & 1242.43   \\
  lc2\_2\_3  & 200           & 6                  & 1944.33       & 2091.35   \\
  lc2\_4\_5s & 400           & 12                 & 4070.48       & 2695.49   \\
  LC2\_8\_5  & 800           & 25                 & 12689.23      & 3820.94   \\
  \hline
\end{tabular}
\end{center}
\end{table}

When we compare the results given by Tables 10 through 13 we can see that the constructive solution gives better results than the random initial solution for both SA and Taboo, but they also take more computational time. In general terms, the best solution as seen in optimal solution is given by the genetic algorithm and Taboo with the constructive solution. The best solution as seen in computational time is given by genetic algorithm and simulated annealing with a random initial solution. 

\section{Conclusions}\label{sec_conclusions}
\begin{itemize}
\item The best solution as seen in optimal solution is given by the genetic algorithm and Taboo with the constructive solution in hybrid methods.
\item The best solution as seen in computational time is given by the greedy approach in constructive methods.
\item An interesting matter with VND is that even though the solution is bad it has a good computational time.
\item The mix between the ALNS and the VND is a kind of hybrid but it is slow and gives bad solution, but the hybrids with the genetic algorithm give better results due to the fact that they use the hybridization to cover its weakness.
\item We can see from all the computational times that when the instances are small, the times between the heuristics are very similar but when the instances are bigger, the times between the heuristics are very different. 
\end{itemize}

\section{Future Work}
An interesting matter that it could be talked about in the future is about the lower bound of this problem, because in none of the articles we have reviewed the topic has risen up, and it's important because we can not compare how much closer we are getting to a better solution, we can see it gets better every time but we can not see if we are closer to our goal. Also, it would be interesting to run the algorithm several times rising the number of transfer points with the purpose of seeing how it affects the computational time and how much it diminishes the objective function.

{\small
\bibliographystyle{authordate1}
\bibliography{main}{\nocite{*}}}

\end{document}